\documentclass[11 pt]{article}
\usepackage{amsmath}
\usepackage{amsfonts}
\usepackage{color}
\begin{document}
\vspace*{.5cm}
\begin{center}
{\large{\bf Screen Transversal Cauchy Riemann Lightlike Submanifolds}}

\bigskip

{\small Bur\c{c}in DO\u{G}AN$^a$, Bayram \c{S}AH\.{I}N$^b$ and Erol YA\c{S}AR$^a$}
\bigskip

{\footnotesize {\it a}:Department of Mathematics, Mersin University, Mersin, Turkey \\
{\it b}: (Corresponding Author)Department of Mathematics, Ege University, \.{I}zmir, Turkey}

{\footnotesize e-mails: bdogan@mersin.edu.tr, bayram.sahin@ymail.com, yerol@mersin.edu.tr}
\end{center}

\vspace{.7cm}

\noindent {\bf Abstract}. {\small
We introduce a new class of lightlike submanifolds, namely, Screen Transversal Cauchy Riemann (STCR)-lightlike submanifolds, of indefinite Kähler manifolds. We show that this new class is an umbrella of screen transversal lightlike, screen transversal totally real lightlike and CR-lightlike submanifolds. We give a few examples of a STCR lightlike submanifold, investigate the integrability of various distributions, obtain a characterization of such lightlike submanifolds in a complex space form and find new conditions for the induced connection to be a metric connection. Moreover, we investigate the existence of totally umbilical (STCR)-lightlike submanifolds and minimal (STCR)-lightlike submanifolds. The paper also contains several examples.}

\section*{\bf 1.~Introduction}
\setcounter{equation}{0}
\renewcommand{\theequation}{1.\arabic{equation}}

In \cite{Duggal-Bejancu}, Duggal and Bejancu studied  the geometry of arbitrary lightlike submanifolds of semi-Riemannian
manifolds. Since then, many authors have studied the geometry of lightlike hypersurfaces
and lightlike submanifolds. Lightlike geometry has its applications in general relativity, particularly in
black hole theory. Indeed, it is known that lightlike hypersurfaces are examples of physical models of
Killing horizons in general relativity \cite{Galloway}. A Killing horizon is a
lightlike hypersurface whose generating null vector can be normalized so as to coincide with one of the
Killing vector. The surface of a black hole is described in terms of Killing horizon. This relation has its
roots in Hawking's area theorem which states that if matter satisfies the dominant energy condition, then
the area of the black hole can not decrease \cite{Hawking}.

On the other hand, complex manifolds, in particular K\"{a}hler manifolds, have been a useful tool in mathematical physics.   Since the 2-form $\rho,$ defined by $\rho(X,Y)=Ric(X,JY),$ on K\"{a}hler
manifold is closed, it represents the first Chern class $C_1.$ Complex manifolds with Ricci flat K\"{a}hler metric are called Calabi-Yau manifolds. The Calabi–
Yau manifolds have their application in super string theory which is based on a 10-dimensional manifold $ M \times V_4,$ where $V_4$ is ordinary spacetime and $M$ is a 6-dimensional manifold which is at least approximately Ricci flat. We also note that, in general, the complex versions of Einstein equations are easier to solve than their real forms \cite{Woodhouse}.

The main  difference between the theory of lightlike submanifolds
and semi-Riemannian submanifolds arises due to the fact that in the
first case, a part of the normal vector bundle $TM^{\perp}$ lies in
the tangent bundle $TM$ of the submanifold $M$ of a semi-Riemannian
manifold $\bar{M},$ whereas in the second case $TM\cap TM^{\perp}
=\{0 \}.$ Thus, the basic problem of lightlike submanifolds is to
replace the intersecting part by a vector subbundle whose sections
are nowhere tangent to $M$.  In \cite{Duggal-Bejancu}, Duggal and Bejancu
introduced a non-degenerate screen distribution to construct a
nonintersecting lightlike transversal vector bundle of the tangent
bundle and then they  studied the geometry of arbitrary lightlike
submanifold of a semi-Riemannian manifold. Although Duggal-Bejancu's approach is extrinsic, there is another approach which is intrinsic and that theory can be found in \cite{Kupeli}. In this paper, we follow Duggal-Bejancu's approach given in \cite{Duggal-Bejancu}.

    In \cite{Duggal-Bejancu}, Duggal and Bejancu defined
CR-lightlike submanifolds of indefinite K\"{a}hler manifolds as a
generalization of lightlike real hypersurfaces of indefinite K\"{a}hler
manifolds. Contrary to the Riemannian CR-submanifolds, CR-lightlike
submanifolds do not contain invariant and totally real lightlike
submanifolds. Therefore, in \cite{Duggal-Sahin-1} Duggal and the second author
defined screen CR-submanifolds of indefinite K\"{a}hler manifolds and
showed that screen CR-submanifolds include invariant submanifolds as
well as screen real submanifolds. Later, in \cite{Duggal-Sahin-2}, the authors gave a
generalization of this notion defining generalized CR-lightlike
submanifolds, obtaining CR-lightlike and
screen CR-lightlike submanifolds as particular cases. Since then, many papers have appeared on the subject, see for instance;  \cite{Gogna-Kumar-Nagaich}, \cite{Gupta-Sharfuddin}, \cite{Gupta-Sharfuddin-Upadhyay}, \cite{Jain-Kumar-Nagaich}, \cite{Kumar-Kumar-Nagaich}, \cite{Sachdeva-Kumar-Bhatia}, \cite{Sachdeva-Kumar-Bhatia2}, \cite{Jain-Kumar-Nagaich2}.

  However, one can observe that  CR-lightlike, screen CR-lightlike and
generalized CR-lightlike  do
not contain real lightlike curves. Therefore, in \cite{Sahin}, the second author  introduced screen transversal
lightlike submanifolds of indefinite K\"{a}hler manifolds and showed that such lightlike submanifolds include real lightlike curves.

In this paper, as a generalization of CR-lightlike submanifolds and screen transversal lightlike submanifolds, we introduce screen transversal CR-lightlike submanifolds and study the geometry of such lightlike submanifolds.
In section 2, we give basic information needed for this paper. In
section~3, we first define {\bf
STCR}-lightlike submanifolds, then  prove a
characterization theorem and investigate the geometry of leaves of distributions which are arisen from definition. In general, the induced connection of a lightlike submanifold is not a metric connection. Therefore it is an important problem to find conditions for the induced connection to be a metric connection. In section 3, we also find necessary and sufficient conditions for the induced connection to be a metric connection. In section~4, we
study totally umbilical proper {\bf STCR}- lightlike submanifolds
and prove some existence theorems. In section~5, we give an example of minimal lightlike submanifolds and obtain certain characterizations. Finally, note that the paper
contains several examples.

\section*{{\bf 2.~ Preliminaries}}

\setcounter{equation}{0}
\renewcommand{\theequation}{2.\arabic{equation}}

~~~~Let
($\overline{M},\bar{g}$) be a $2k$-dimensional semi-Riemannian manifold with the semi-Riemannian metric $\bar{g}.$  Denote the constant index of $\bar{g}$ by $q,$ $0<q<2k.$ A tensor field $\bar{J}$ of type (1,1) on $\bar{M}$ is an almost complex structure on $\bar{M}$ if $\bar{J}^2=-I,$ $\forall p\in \bar{M},$ where $I$ denotes the identity transformation of $T_{p} \bar{M}$ and such a manifold $\bar{M}$ is called an almost complex manifold. Let $\bar{g},$ be a semi-Riemannian metric on an almost complex manifold $\bar{M}$ such that
\begin{equation}
\bar{g}(X,Y)=\bar{g}(\bar{J}X,\bar{J}Y) \label{eq:2.1} \,\,\,,\,\ \forall
X,Y \in \Gamma(T\bar{M})
\end{equation}
is satisfied. Then, $\bar{g}$ is a Hermitian metric and $(\bar{M},\bar{J},\bar{g})$ is an indefinite almost Hermitian manifold. Denote the Levi-Civita connection on an indefinite almost Hermitian manifold $\bar{M}$ with respect to $\bar{g}$ by $\bar{\nabla}.$ If $\bar{J}$ is parallel with respect to $\bar{\nabla}$, i.e.,
\begin{equation}
(\bar{\nabla}_{X} \bar{J})Y=0 \label{eq:2.2} \,\,\,,\,\ \forall X,Y \in\Gamma(T\bar{M}),
\end{equation}
then $\bar{M}$ is called an indefinite K\"{a}hler manifold. \\
An indefinite complex space form $\bar{M}(c)$ is a connected indefinite K\"{a}hler manifold of constant holomorphic sectional curvature $c$ and its curvature tensor field is calculated as
\begin{eqnarray}
  \bar{R}(X,Y)Z &=&\frac{c}{4}\{\bar{g}(Y,Z)X-\bar{g}(X,Z)Y+\bar{g}(\bar{J}Y,Z)\bar{J}X-\bar{g}(\bar{J}X,Z)\bar{J}Y\nonumber\\
  &+&2\bar{g}(X,\bar{J}Y)\bar{J}Z\}, \label{eq:2.3}
\end{eqnarray}
for any $X,Y\in\Gamma(T\bar{M})$ \cite{Barros-Romero}. \\

~~~~From now on, we use the same notations and formulas in \cite{Duggal-Bejancu}. \\

~~~~Let $(\bar{M},\bar{g})$ be a $(m+n)$-dimensional semi-Riemannian manifold and $(M,g)$ be a $m$-dimensional submanifold of $(\bar{M},\bar{g}).$ The induced metric $g$ on $M$ from $\bar{g}$ on $\bar{M}$ does not always have to be non-degenerate. If the induced metric $g$ is degenerate on $M$ and $rank(Rad(TM))=r,$ $1 \leq r \leq m,$ then $(M,g)$ is called a lightlike submanifold of $(\bar{M},\bar{g}),$ where the {\sl radical distribution} $Rad(TM)$ and the {\sl normal bundle} $TM^{\perp}$ of the tangent bundle $TM$ are defined as $$Rad(TM) = TM \cap TM^{\perp}$$ and $$TM^{\perp} = \cup_{x \in M} \{ u \in T_x\,\bar{M} \mid \bar{g}(u,v) = 0,\,\,\ \forall v \in  T_x\,M\}.$$
Since $TM$ and $TM^{\perp}$ are degenerate vector subbundles, there exist complementary non-degenerate distributions $S(TM)$ and $S(TM^{\bot})$ of $Rad(TM)$ in $TM$ and $TM^{\perp},$ respectively, which are called the {\sl screen distribution} and {\sl screen transversal bundle (or co-screen distribution)} of $M$ such that
$$TM = S(TM) \perp Rad(TM)\,,\,\,\  TM^{\bot}=S(TM^{\bot}) \perp Rad(TM).$$
On the other hand, consider a orthogonal complementary bundle $S(TM)^{\bot}$ to $S(TM)$ in $T\bar{M}$ such that
\begin{equation}
S(TM)^{\bot}=S(TM^{\bot})\bot S(TM^{\bot})^{\bot} \label{eq:2.a}
\end{equation}
where $S(TM^{\bot})^{\bot}$ is the orthogonal complementary to $S(TM^{\bot})$ in $S(TM)^{\bot}.$ We now recall the following important result. \\

\noindent{\bf Theorem~2.1.~}{\sl Let $(M,g,S(TM),S(TM^{\bot}))$ be a $r$-lightlike submanifold of a semi-Riemannian manifold $(\bar{M},\bar{g}).$ Then, there exists a complementary vector bundle $ltr(TM)$ called a lightlike transversal bundle of $Rad(TM)$ in $S(TM^{\bot})^{\bot}$ and a basis of $\Gamma(ltr(TM)\mid_{U})$ consists of smooth sections $\{N_1,...,N_r\}$ of $S(TM^{\bot})^{\bot}\mid_{U}$ such that
$$\bar{g}(\xi_i,N_j) = \delta_{ij}\,,\,\,\ \bar{g}(N_i, N_j)=0\,,\,\,\ i,j=1,..,r, $$
where $\{\xi_1,...,\xi_r\}$ is a basis of $\Gamma(Rad(TM))$} \cite[page 144]{Duggal-Bejancu}.\\

~~~~This result implies that there exists a complementary (but not orthogonal) vector bundle $tr(TM)$ to $TM$ in $T\bar{M}|_M,$ which called {\sl transversal vector bundle}, such that the following decompositions are hold:
\begin{equation}
tr(TM)= ltr(TM) \perp S(TM^{\perp}) \label{eq:2.b}
\end{equation}
and
\begin{equation}
S(TM^{\bot})^{\bot}=Rad(TM)\oplus ltr(TM). \label{eq:2.c}
\end{equation}
Thus, using (\ref{eq:2.a}), (\ref{eq:2.b}) and (\ref{eq:2.c}) we get
\begin{align}
T\bar{M}|_M &=S(TM) \perp S(TM)^{\bot} \nonumber \\
&=S(TM) \perp [Rad(TM) \oplus ltr(TM)] \perp S(TM^{\perp}) \nonumber \\
&=TM \oplus tr(TM). \label{eq:2.d}
\end{align}

A submanifold $(M, g, S(TM), S(TM^{\perp})$ is called \\
(1):\, $r$\,-\,lightlike if $r\,< min\{m,\,n\},$ \\
(2):\, Co\,-\,isotropic if $r\,=\,n\,<\,m,\,\,\ i.e.,\,\,\ S(TM^{\perp}) = \{0\},$\\
(3):\, Isotropic if $r\,=\,m\,<\,n,\,\,\ i.e.,\,\,\ S(TM) = \{0\}$ and \\
(4):\, Totally lightlike if $r\,=\,m\,=\,n,\,\,\ i.e.,\,\,\ S(TM) = \{0\}= S(TM^{\perp}).$  \\

The Gauss and Weingarten equations of $M$ are given by
\begin{equation}
\bar{\nabla}_{X}Y=\nabla_{X}Y+h(X,Y)\,\,,\,\,\ \forall
X,Y\in\Gamma(TM) \label{eq:2.4}
\end{equation}
and
\begin{equation}
\bar{\nabla}_{X}V=-A_{V}X+\nabla^{t}_{X}V \,\,\,,\,\,\ \forall
X\in\Gamma(TM)\,\,,\,\ V\in\Gamma(tr(TM)), \label{eq:2.5}
\end{equation}
where $\{\nabla_{X}Y,A_{V}X\}$ and $\{ h(X,Y),\nabla^{t}_{X}V\}$ are belong to $\Gamma(TM)$ and $\Gamma(tr(TM)),$ respectively. $\nabla$ and $\nabla^t$ are linear connections on $M$ and on the vector bundle $tr(TM),$ respectively. The second fundamental form $h$ is a symmetric ${\cal F}(M)$-bilinear form on $\Gamma(TM)$ with values in $\Gamma(tr(TM))$ and the shape operator $A_V$ is a linear endomorphism of $\Gamma(TM).$ If we consider (\ref{eq:2.d}) and using the projectors
$$L:tr(TM)\rightarrow ltr(TM)\,,\,\,\,\,\ S:tr(TM)\rightarrow S(TM^\perp),$$ we can write, for $\forall X,Y\in\Gamma(TM),$ $N\in\Gamma(ltr(TM))$ and $W\in\Gamma(S(TM^\perp)),$
\begin{eqnarray}
\bar{\nabla}_{X}Y &=& \nabla_{X}Y+h^{\l}(X,Y)+h^s(X,Y), \label{eq:2.6}\\
\bar{\nabla}_{X}N &=& -A_{N}X+\nabla^{\l}_{X}(N)+D^{s}(X,N), \label{eq:2.7}\\
\bar{\nabla}_{X}W &=& -A_{W}X+D^{l}(X,W)+\nabla^{s}_{X}(W), \label{eq:2.8}
\end{eqnarray}
where $\{\nabla^{\l}_{X}N,D^{l}(X,W)\},$ $\{D^{s}(X,N),\nabla^{s}_{X}W\}$ are parts of $ltr(TM),$ $S(TM^{\bot}),$ respectively and $h^{\l}(X,Y)=Lh(X,Y)\in\Gamma(ltr(TM)),$ $h^s(X,Y)=Sh(X,Y)\in\Gamma(S(TM^{\bot}))$. Denote the projection of $TM$ on $S(TM)$ by $\bar{P}.$ Then, by using (\ref{eq:2.4}), (\ref{eq:2.6})-(\ref{eq:2.8}) and taking account
that $\bar{\nabla}$ is a metric connection we obtain
\begin{equation}
\bar{g}(h^s(X,Y),W)+\bar{g}(Y,D^l(X,W))=g(A_{W}X,Y), \label{eq:2.9}
\end{equation}
\begin{equation}
\bar{g}(D^s(X,N),W)=\bar{g}(N,A_{W}X), \label{eq:2.10}
\end{equation}
\begin{equation}
{\nabla}_{X}{\bar{P}}Y = {\nabla}^*_{X}{\bar{P}}Y+h^*(X,{\bar{P}}Y)
\label{eq:2.11}
\end{equation}
and
\begin{equation}
{\nabla}_{X}{\xi} = -A^*_{\xi}X+{{\nabla}^*}^t_{X}{\xi}
\label{eq:2.12}
\end{equation}
for $X,Y\in \Gamma(TM)$ and $\xi \in \Gamma(RadTM),$ where ${\nabla}^*$ and ${{\nabla}^*}^t$ are induced connections on $S(TM)$ and $Rad(TM)$. On the other hand, $h^{*}$ and $A^{*}$ are $\Gamma(Rad(TM))$-valued and $\Gamma(S(TM))$-valued ${\cal F}(M)$-bilinear forms on $\Gamma(TM)\times \Gamma(S(TM))$ and $\Gamma(Rad(TM))\times \Gamma(TM),$ respectively. $h^{*}$ is called lokal second fundamental form on $S(TM)$ and $A^{*}$ is second fundamental form of $Rad(TM).$ By using above equations we obtain
\begin{eqnarray}
\bar{g}(h^l(X,{\bar{P}}Y),\xi)&=&g(A^*_{\xi}X,{\bar{P}}Y), \label{eq:2.13} \\
\bar{g}(h^*(X,{\bar{P}}Y),N)&=&g(A_{N}X,{\bar{P}}Y), \label{eq:2.14}\\
\bar{g}(h^l(X,\xi),\xi)=0,& &A^*_{\xi}{\xi}=0. \label{eq:2.15}
\end{eqnarray}

 In general, the induced connection $\nabla$ on $M$ is not metric connection. Since $\bar{\nabla}$ is a metric connection, by using (\ref{eq:2.6}) we get
\begin{equation}
(\nabla_{X}g)(Y,Z)=\bar{g}(h^l(X,Y),Z)+\bar{g}(h^l(X,Z),Y).
\label{eq:2.16}
\end{equation}
However, it is important to note that $\nabla^{\star}$ is a metric connection on $S(TM).$ We denote curvature tensor of a lightlike submanifold by $R,$ then the Gauss equation for lightlike submanifolds is given by
\begin{eqnarray}
  \bar{R}(X,Y)Z&=&R(X,Y)Z+A_{h^l(X,Z)}Y-A_{h^l(Y,Z)}+A_{h^s(X,Z)}Y \nonumber\\
  &-& A_{h^s(Y,Z)}X+(\nabla_{X}h^l)(Y,Z)-(\nabla_{Y}h^l)(X,Z) \nonumber\\
  &+&D^l(X,h^s(Y,Z))-D^l(Y,h^s(X,Z)) \nonumber\\
  &+&(\nabla_{X}h^s)(Y,Z)-(\nabla_{Y}h^s)(X,Z) \nonumber\\
  &+&D^s(X,h^l(Y,Z))-D^s(Y,h^l(X,Z)), \label{eq:2.17}
\end{eqnarray}
for any $X,Y,Z\in\Gamma(TM).$

\section*{{\bf 3.\, Screen Transversal Cauchy Riemann Lightlike Submanifolds}}

\setcounter{equation}{0}
\renewcommand{\theequation}{3.\arabic{equation}}

As we have mentioned in the introduction, CR-lightlike submanifolds and generalized CR-lightlike submanifolds of an indefinite K\"{a}hler manifold includes real lightlike hypersurfaces, however such lightlike submanifolds excludes real lightlike curves. On the other hand, screen transversal lightlike submanifolds covers real lightlike curves, however this class does not include real lightlike hypersurfaces. But real lightlike curves and real lightlike hypersurfaces are important subjects in relativity theory. Indeed, the propagation of light and other zero rest mass particles are described by the null geodesics of the spacetime \cite{Flaherty} and lightlike hypersurfaces are examples of physical models of Killing horizons in general relativity. Therefore we ask the following question \\

{\it Are there any lighlike submanifolds of indefinite K\"{a}hler manifolds containing both real lightlike curves and real lightlike hypersurfaces?}\\

To give affirmative answer to above question, in this section, we introduce screen transversal Cauchy Riemann lightlike submanifolds of an indefinite K\"{a}hler manifold as a generalization of CR-lightlike submanifolds and screen transversal lightlike submanifolds. We give examples, obtain a characterization and find necessary and sufficient conditions for the induced connection, which is not metric connection in general, to be a metric connection. We also check the effect of the notion of mixed geodesic on the geometry of submanifolds.\\

\noindent{\bf Definition~1.~}{\sl Let $M$ be a real $r$-lightlike submanifold of an indefinite K\"{a}hler manifold $\bar{M}.$ Then we say that $M$ is a screen transversal Cauchy Riemann (STCR) lightlike submanifold if the following conditions are satisfied:}
\begin{enumerate}
  \item [(A)] {\sl There exist two subbundles $D_1$ and $D_2$ of
  $Rad(TM)$ such that }
  \begin{equation}
  Rad(TM)=D_1 \oplus D_2 \,,\,\ \bar{J}(D_1)\subset S(TM)\,,\,\ \bar{J}(D_2) \subset
  S(TM^{\perp}), \label{eq:3.1}
  \end{equation}
  \item [(B)] {\sl There exist two subbundles $D_0$ and $D'$ of
  $S(TM)$ such that }
  \begin{equation}
  S(TM)=\{\bar{J}D_1 \oplus D'\} \perp
  D_0 \,,\,\,\ \bar{J}(D_0)=D_0 \,,\,\,\ \bar{J}(D')=L_1\perp S,
  \label{eq:3.2}
  \end{equation}
  {\sl where $D_0$ is a non-degenerate distribution on $M,$ $L_1$ and $S$ are vector subbundles of $ltr(TM)$ and $S(TM^{\perp}),$ respectively.}
\end{enumerate}

From definition of a screen transversal Cauchy Riemann lightlike submanifold, we obtain that the tangent bundle of a screen transversal Cauchy Riemann lightlike submanifold is decomposed as follows
\begin{equation}
TM=D \oplus \tilde{D}, \label{eq:3.3}
\end{equation}
where
\begin{equation}
D=D_0 \oplus D_1 \oplus \bar{J}D_1\ \label{eq:3.4}
\end{equation}
and \
\begin{equation}
\tilde{D}=D_2 \oplus \bar{J}L_1 \oplus \bar{J}S. \label{eq:3.5}
\end{equation}
It is clear that $D$ is invariant and $\tilde{D}$ is anti-invariant. Furthermore, we have
\begin{equation*}
ltr(TM)=L_1\oplus L_2 \,,\,\ \bar{J}(L_1)\subset S(TM)\,,\,\ \bar{J}(L_2) \subset S(TM^{\perp}),
\end{equation*}
and
\begin{equation*}
S(TM^{\bot})=\{\bar{J}(D_2)\oplus \bar{J}(L_2)\}\bot S.
\end{equation*}
We say that $M$ is a proper screen transversal Cauchy Riemann lightlike submanifold
of an indefinite K\"{a}hler manifold if $D_1 \neq \{0\}\,,\,\,\ D_2 \neq\{0\} \,,\,\,\ D_0 \neq \{0\}$ and $S \neq\{0\}.$ For proper  screen transversal Cauchy Riemann lightlike submanifold we note that the following features:
\begin{enumerate}
  \item The condition (A) implies that $dim(Rad(TM)) \geq 2.$
  \item The  condition (B) implies $dim(D)=2s \geq 4,$ $dim(D')\geq
  2$ and $dim(D_2)=dim(L_2).$  Thus $dim(M) \geq 7$ and
  $dim(\bar{M})\geq 12.$
  \item Any proper $7-$ dimensional screen transversal Cauchy Riemann lightlike
  submanifold must be $2-$lightlike.
  \item (A) and K\"{a}hler manifold $\bar{M}$ imply that $index(\bar{M}) \geq 2.$
\end{enumerate}

\noindent{\bf Proposition~3.1.~}{\sl A STCR lightlike submanifold $M$ of an indefinite K\"{a}hler manifold $\bar{M}$ is a CR-lightlike submanifold (respectively, screen transversal lightlike submanifold) if and only if $D_2=\{0\}$ (respectively, $D_1=\{0\}.$)}\\
\noindent{\bf Proof.~} Let $M$ be a CR-lightlike submanifold of an indefinite K\"{a}hler manifold. Then $\bar{J}(Rad(TM))$ is a distribution on $M$ such that $\bar{J}(Rad(TM)) \cap Rad(TM)=\{0\}.$ Thus we obtain $D_1=Rad(TM)$ and $D_2=\{0\}.$ Hence we conclude that $\bar{J}(ltr(TM)) \cap ltr(TM)=\{0\}.$ Then it follows that $\bar{J}(ltr(TM))\subset S(TM).$ Conversely, suppose that $M$ be a STCR lightlike submanifold such that $D_2=\{0\}.$ Then we have $D_1=Rad(TM).$ Hence $\bar{J}(Rad(TM)) \cap Rad(TM)=\{0\},$ that is $\bar{J}(Rad(TM))$ is a  vector subbundle of $S(TM).$ Thus $M$ is a CR-lightlike submanifold. The other assertion can be proved in a similar way. \\

\noindent{\bf Example~1.~} {\sl Let  $M$ be a submanifold of $R^{12}_4$ given by equations
\begin{align*}
x_1&=\sin{u_{2}}\,,\,\,\ x_2=-\cos{u_{2}}\,,\,\,\ x_3=u_{1}\,,\,\,\ x_4=u_{3}-\frac {u_4} {2}\,,\,\,\ x_5=u_{2}, \\ x_6&=0\,,\,\,\ x_7=u_{1}\,,\,\,\ x_8=u_{3}+\frac{u_4}{2}\,,\,\,\ x_9=u_{5}+u_{7}\,,\,\,\ x_{10}=u_{6}-u_{7}, \\ x_{11}&=u_{5}-u_{7}\,,\,\,\ x_{12}=u_{6}+u_{7}.
\end{align*}
Then $TM$ is spanned by
$\{Z_1,Z_2,Z_3,Z_4,Z_5,Z_6,Z_7\}$ where
\begin{align*}
Z_1&={\partial\,x_3}+{\partial\,x_7} \,,\,\,\,\
Z_2=\cos{u_{2}}{\partial\,x_1}+\sin{u_{2}}{\partial\,x_2}+{\partial\,x_5},\\
Z_3&={\partial\,x_4}+{\partial\,x_8}\,,\,\,\,\ Z_4=\frac{1}{2}\{-{\partial\,x_4}+{\partial\,x_8}\},\\
Z_5&={\partial\,x_9}+{\partial\,x_{11}}\,,\,\,\,\ Z_6={\partial\,x_{10}}+{\partial\,x_{12}},\\
Z_7&={\partial\,x_9}-{\partial\,x_{10}}-{\partial\,x_{11}}+{\partial\,x_{12}}.
\end{align*}
Hence $M$ is a $2-$ lightlike submanifold of $R^{12}_4$ with $Rad(TM)=Span\{Z_1,Z_2\}.$ It is easy to see $\bar{J}Z_1=Z_3\in \Gamma(S(TM)),$ thus $D_1=Span\{Z_1\}$ and $D_2=Span\{Z_2\}.$ On the other hand, since $\bar{J}Z_5=Z_6 \in \Gamma(S(TM)),$ we obtain $D_0=Span\{Z_5,Z_6\}$ and by direct calculations, we get the lightlike transversal bundle spanned by
$$N_1=\frac{1}{2}\{-{\partial\,x_3}+{\partial\,x_7}\}\,,\,\,\
N_2=\frac{1}{2}\{-\cos{u_{2}}{\partial\,x_1}-\sin{u_{2}}{\partial\,x_2}+{\partial\,x_5}\}.$$
Then we see that $L_1=Span\{N_1\}$, $L_2=Span\{N_2\}$, $S(TM^{\perp})=Span\{\bar{J}Z_2,\bar{J}N_2,\bar{J}Z_7\}$ and $S=Span\{\bar{J}Z_7=W\}.$ Thus, $D'=Span\{\bar{J}N_1=Z_4,\bar{J}Z_7=W\}$ and $M$ is a proper STCR lightlike submanifold. }\\

\noindent{\bf Example~2.~} {\sl Every CR-lightlike submanifold of an indefinite K\"{a}hler manifold is a STCR lightlike submanifold with $D_2=\{0\}.$}\\

It is known that every real lightlike hypersurface is a CR-lightlike submanifold \cite{Duggal-Bejancu}, therefore a real lightlike hypersurface is an example of STCR lightlike submanifold.\\

\noindent{\bf Example~3.~} {\sl Every screen transversal lightlike submanifold of an indefinite K\"{a}hler manifold is a STCR lightlike submanifold with $D_1=\{0\}.$}\\

It is known that every real lightlike curve is an isotropic screen transversal lightlike submaniold \cite{Sahin}, therefore a real lightlike curve is an example of STCR lightlike submanifold.\\

\noindent{\bf Proposition~3.2.~}{\sl There exist no coisotropic, isotropic or totally lightlike proper STCR lightlike submanifolds $M$ of an indefinite K\"{a}hler manifold. Any isotropic STCR lightlike submanifold is a screen transversal lightlike submanifold. Also, a coisotropic STCR lightlike submanifold is a CR-lightlike submanifold.}\\
\noindent{\bf Proof.~} Let $M$ be a proper STCR lightlike submanifold. From definition of proper STCR lightlike submanifold, we know that $D_1 \neq \{0\}\,,\,\ D_2 \neq\{0\} \,,\,\ D_0 \neq \{0\}$ and $S \neq\{0\},$ i.e., both $S(TM)$ and $S(TM^{\bot})$ are non-zero. Thus, $M$ can not be a coisotropic, isotropic or totally lightlike submanifold. On the other hand, if $M$ be a isotropic STCR lightlike submanifold, then $S(TM)=\{0\},$ i.e., $\bar{J}D_1=\{0\}$ and $Rad(TM)=D_2.$ Thus, we get $\bar{J}Rad(TM)=\bar{J}(D_2)\subset \Gamma(S(TM^{\bot}))$ and $M$ is a screen transversal lightlike submanifold. Similarly, if $M$ is a coisotropic STCR lightlike submanifold, then $S(TM^{\bot})=\{0\},$ i.e., $\bar{J}D_2=\{0\}$ and $Rad(TM)=D_1.$ Since, $\bar{J}Rad(TM)=\bar{J}(D_1)\subset \Gamma(S(TM))$ then $M$ is a CR-lightlike submanifold. \\

Now, we denote the projections from $\Gamma (TM)$ to $\Gamma(D_0),$ $\Gamma(\bar{J}D_1),$ $\Gamma(\bar{J}L_1),$ $\Gamma(\bar{J}S),$ $\Gamma(D_1)$ and $\Gamma(D_2)$ by $P_0,$ $ P_1,$ $P_2,$ $P_3,$ $S_1$ and $ S_2,$ respectively. We also denote the projections from $\Gamma(tr(TM))$ to $\Gamma(\bar{J}D_2),$ $\Gamma(\bar{J}L_2),$ $\Gamma(S),$ $\Gamma(L_1)$ and $\Gamma(L_2)$ by $R_1,$ $R_2,$ $R_3,$ $Q_1$ and $Q_2,$ respectively.
Thus, we write
\begin{align}
X&=PX+QX =P_0 X+P_1 X+P_2 X+P_3 X+S_1 X+S_2 X \label{eq:3.6}
\end{align}
and
\begin{equation}
\bar{J}X=TX+\omega X, \label{eq:3.7}
\end{equation}
for $X \in \Gamma(TM),$
where $PX\in\Gamma(D),$ $QX\in\Gamma(\tilde{D})$ and $TX$ and $wX$ are the tangential parts and the transversal parts of $\bar{J}X,$ respectively.
Applying $\bar{J}$ to (\ref{eq:3.6}) and denoting $\bar{J}P_0,$ $\bar{J}P_1,$ $\bar{J}P_2,$ $\bar{J}P_3,$ $\bar{J}S_1,$ $\bar{J}S_2$ by $T_0,$ $T_1,$ $\omega_L,$ $\omega_S,$ $T_{\bar{1}},$ $\omega_{\bar{2}},$ respectively, we obtain
\begin{eqnarray}
\bar{J}X&=&T_0 X+T_1 X+T_{\bar{1}} X+\omega_L X+\omega_S X+\omega_{\bar{2}} X, \label{eq:3.8}
\end{eqnarray}
for $X \in \Gamma(TM),$ where $T_0 X \in \Gamma(D_0),$ $T_1 X \in \Gamma(D_1),$ $T_{\bar{1}} X \in \Gamma(\bar{J}D_1),$ $\omega_L X \in \Gamma(L_1),$ $\omega_S X \in \Gamma(S),$ and $\omega_{\bar{2}} X \in \Gamma(\bar{J}D_2).$ Similarly we can write, for any $V \in \Gamma(tr(TM)),$
\begin{equation}
V=R_1 V+R_2 V+R_3 V+Q_1 V+Q_2 V \label{eq:3.9}
\end{equation}
and we denote $\bar{J}R_1,$ $\bar{J}R_2,$ $\bar{J}R_3,$ $\bar{J}Q_1,$ $\bar{J}Q_2$ by $B_2,$ $C_L,$ $B_{\bar{S}},$ $B_{\bar{L}},$ $C_{\bar{L}},$ respectively, we write
\begin{align}
\bar{J}V&=B_2 V+B_{\bar{S}} V+B_{\bar{L}} V+C_L V+C_{\bar{L}} V, \label{eq:3.10}
\end{align}
where $BV$ and $CV$ are sections of $TM$ and $tr(TM),$ respectively. Now, differentiating (\ref{eq:3.8}) and using (\ref{eq:2.2}), (\ref{eq:2.4}), (\ref{eq:2.6})-(\ref{eq:2.8}) and (\ref{eq:3.10}), $\forall X,Y\in\Gamma(TM),$ we have
\begin{align*}
&\nabla_X TY+h^{l}(X,TY)+h^{s}(X,TY)+\{-A_{\omega_L Y} X+\nabla_{X}^{l}(\omega_L Y)+D^{s}(X,\omega_L Y)\} \\
&\,\,\,\,\,\,\,\,\,\,\,\,\,\,\,\,\,\,\,\,\,\,\,\,\,\,\,\,\,\,\,\,\,\,\,\,\,\,\,\,\,\,\,\,\,\,\,\,\,\,\,\,\,\,\,\,\,\,\,\,\,\,\,\,\,\,\,\,\,\,\,\,\,\,\,\,\,\,\,\,\,\,\,\,\,\ +\{-A_{\omega_S Y} X+\nabla_{X}^{s}(\omega_S Y)+D^{l}(X,\omega_S Y)\} \\
&\,\,\,\,\,\,\,\,\,\,\,\,\,\,\,\,\,\,\,\,\,\,\,\,\,\,\,\,\,\,\,\,\,\,\,\,\,\,\,\,\,\,\,\,\,\,\,\,\,\,\,\,\,\,\,\,\,\,\,\,\,\,\,\,\,\,\,\,\,\,\,\,\,\,\,\,\,\,\,\,\,\,\,\,\,\ +\{-A_{\omega_{\bar{2}} Y} X+\nabla_{X}^{s}({\omega_{\bar{2}} Y})+D^{l}(X,{\omega_{\bar{2}} Y})\}\\
&=T\nabla_X Y +\omega_L\nabla_X Y+\omega_S\nabla_X Y+\omega_{\bar{2}}\nabla_X Y+Bh^{l}(X,Y)+Ch^{l}(X,Y) \\
&\,\,\,\,\,\,\,\,\,\,\,\,\,\,\,\,\,\,\,\,\,\,\,\,\,\,\,\,\,\,\,\,\,\,\,\,\,\,\,\,\,\,\,\,\,\,\,\,\,\,\,\,\,\,\,\,\,\,\,\,\,\,\,\,\,\,\,\,\,\,\,\,\,\,\,\,\,\,\,\,\,\,\,\,\,\,\,\,\,\,\,\,\ +Bh^{s}(X,Y)+Ch^{s}(X,Y).
\end{align*}
Considering the tangential, lightlike transversal and screen transversal parts of this equation we obtain
\begin{align}
\nabla_X TY - T\nabla_X Y &=(\nabla_X T)Y \nonumber \\
&=A_{\omega_L Y}X+A_{\omega_S Y}X+A_{\omega_{\bar{2}} Y}X+Bh(X,Y), \label{eq:3.11}
\end{align}
\begin{align}
D^{l}(X,\omega_S Y)+D^{l}(X,\omega_{\bar{2}} Y)&=\omega_L(\nabla_X Y)-\nabla_X ^{l}(\omega_L Y) \nonumber \\
&-h^{l}(X,TY)+Ch^{l}(X,Y) \label{eq:3.12}
\end{align}
and
\begin{align}
D^{s}(X,\omega_L Y)&=\omega_S(\nabla_X Y)+\omega_{\bar{2}}(\nabla_X Y)-\nabla_X ^{s}(\omega_S Y)\nonumber \\
&-\nabla_X ^{s}(\omega_{\bar{2}} Y)-h^{s}(X,TY)+Ch^{s}(X,Y), \label{eq:3.13}
\end{align}
respectively.\\

The following theorem gives new conditions for the induced connection to be metric connection.\\

\noindent{\bf Theorem~3.1.~}{\sl Let $M$ be a STCR lightlike submanifold of an indefinite Kaehler manifold $\bar{M}.$ Then the induced connection is a metric connection if and only if, $\forall X\in\Gamma(TM),$ the followings are hold:
$$
\nabla^*_X \bar{J}Y+h^*(X,\bar{J}Y) \in \Gamma(\bar{J}D_1) \,,\,\,\ Bh(X,\bar{J}Y)=0 \,,\,\,\,\ Y \in \Gamma(D_1),
$$
$$
A_{\bar{J}Y}X \in \Gamma(\bar{J}D_1) \,,\,\,\ B(\nabla^{s}_X \bar{J}Y+D^{l}(X,\bar{J}Y))=0 \,,\,\,\,\ Y \in \Gamma(D_2).
$$ }
\noindent{\bf Proof.~} For $Y \in \Gamma(RadTM)$ and $ X \in \Gamma(TM),$ from (\ref{eq:2.2}) we can write
$$\bar{\nabla}_X Y=-\bar{J}(\bar{\nabla}_X \bar{J}Y)$$
and using (\ref{eq:2.6}) we get
\begin{equation}
\nabla_X Y+h(X,Y)=-\bar{J}(\nabla_X \bar{J}Y+h(X,\bar{J}Y)).
\label{eq:3.14}
\end{equation}
Since $RadTM=D_1 \oplus D_2,$ for $Y \in \Gamma(D_1),$ from (\ref{eq:2.11}), (\ref{eq:3.7}) and (\ref{eq:3.10}) we have
\begin{equation*}
\nabla_X Y+h(X,Y)=-T(\nabla_X \bar{J}Y)-\omega(\nabla_X \bar{J}Y)-B(h(X,\bar{J}Y))-C(h(X,\bar{J}Y))
\end{equation*}
and from (\ref{eq:2.11})
\begin{eqnarray*}
\nabla_X Y+h(X,Y)&=&-T\nabla^{*}_X \bar{J}Y-w\nabla^{*}_X \bar{J}Y-Th^{*}(X,\bar{J}Y)-wh^{*}(X,\bar{J}Y) \\
&-&Bh(X,\bar{J}Y)-Ch(X,\bar{J}Y)
\end{eqnarray*}
is obtained. Taking the tangential parts of this equation we derive
\begin{align}
\nabla_X Y&=-T(\nabla^{*}_X \bar{J}Y+h^{*}(X,\bar{J}Y))-Bh(X,\bar{J}Y). \label{eq:3.15}
\end{align}
In similar way, for $X \in \Gamma(TM)$ and $Y \in \Gamma(D_2),$ using (\ref{eq:2.8}) we get
\begin{align}
\nabla_{X} Y&=TA_{\bar{J}Y} X-B(\nabla^{s}_{X} \bar{J}Y+D^{l}(X,\bar{J}Y)). \label{eq:3.16}
\end{align}
Then our assertion follows from (\ref{eq:3.15}), (\ref{eq:3.16}) and Theorem~2.4 in [4, p.161].  \\

\noindent{\bf REMARK~1.~}{\sl It follows from Theorem~5.1 of \cite[page 49]{Duggal-Bejancu} that, under the conditions of Theorem~3.1, $Rad(TM)$ of this class of STCR lightlike submanifolds is an integrable Killing distribution.} \\

Blair and Chen \cite{Blair-Chen} have obtained a characterization of Riemannian CR submanifolds of a complex space form $\bar{M}(c)$ with $c\neq0$.  Here we give a characterization of STCR lightlike submanifolds in an indefinite complex space form in terms of the curvature tensor field of the ambient space. \\

\noindent{\bf Theorem~3.2.~}{\sl A lightlike submanifold $M$ of an indefinite complex space form $\bar{M}(c)$ with $c\neq0$ is a STCR lightlike submanifold with $D_0\neq0,$ iff
\begin{enumerate}
\item [(a)] {\sl The maximal complex subspaces of $T_{p}M \,,\,\ p\in M,$ define a distribution}
\begin{equation*}
D=D_0\bot D_1\bot \bar{J}(D_1),
\end{equation*}
where $Rad(TM)=D_1\bot D_2$ and $D_0$ is a non-degenerate complex distribution.
\item [(b)] {\sl There exists a lightlike transversal vector bundle $ltr(TM)$ such that for $\forall X,Y\in\Gamma(D)\,,\,\ N_{1}, N_{2}\in\Gamma(ltr(TM)),$}
\begin{equation*}
\bar{g}(\bar{R}(X,Y)N_{1},N_{2})=0.
\end{equation*}
\item [(c)] {\sl There exists a vector subbundle $\bar{J}(S)$ on $M$ such that for $\forall X,Y\in\Gamma(D)\,,\,\ W_{1}, W_{2}\in\Gamma(\bar{J}(S)),$}
\begin{equation*}
\bar{g}(\bar{R}(X,Y)W_{1},W_{2})=0,
\end{equation*}
 where $\bar{J}(S)$ is orthogonal to $D$ and $\bar{R}$ is the curvature tensor of $\bar{M}(c).$
\end{enumerate} }
\noindent{\bf Proof.~} $\Rightarrow):$ Let $M$ be a STCR lightlike submanifold of an indefinite K\"{a}hler manifold $\bar{M}.$ Since
\begin{equation*}
 D=D_0\bot D_1\bot \bar{J}(D_1)
\end{equation*}
(a) holds. On the other hand, for $\forall X,Y\in\Gamma(D)$ and $N_1\,,\ N_2\in\Gamma(ltr(TM)),$ from (\ref{eq:2.3}) we get
\begin{align*}
\bar{g}(\bar{R}(X,Y)N_1,N_2)&=\frac{c}{4}\{\bar{g}(Y,N_1)\bar{g}(X,N_2)
-\bar{g}(X,N_1)\bar{g}(Y,N_2) \\
&+\bar{g}(\bar{J}Y,N_1)\bar{g}(\bar{J}X,N_2)
-\bar{g}(\bar{J}X,N_1)\bar{g}(\bar{J}Y,N_2) \\
&+2\bar{g}(X,\bar{J}Y)\bar{g}
(\bar{J}N_1,N_2)\} \\
&=0.
\end{align*}
Thus (b) holds. In a similar way, since for $\forall X,Y\in\Gamma(D)$ and $W_1\,,\ W_2\in\Gamma(\bar{J}(S)),$ we have
\begin{align*}
\bar{g}(\bar{R}(X,Y)W_1,W_2)&=\frac{c}{4}\{\bar{g}(Y,W_1)\bar{g}(X,W_2)
-\bar{g}(X,W_1)\bar{g}(Y,W_2) \\
&+\bar{g}(\bar{J}Y,W_1)\bar{g}(\bar{J}X,W_2)
-\bar{g}(\bar{J}X,W_1)\bar{g}(\bar{J}Y,W_2) \\
&+2\bar{g}(X,\bar{J}Y)\bar{g}
(\bar{J}W_1,W_2)\} \\
&=0.
\end{align*}
Then (c) also holds. \\
$\Leftarrow):$ Conversely, we assume that (a), (b) and (c) are provided. From (a), for the maximal distribution $D=D_0 \bot D_1 \bot \bar{J}(D_1)$ and $Rad(TM)=D_1 \bot D_2,$ while $\bar{J}(D_1)$ is a invariant distribution on $TM,$ $D_2$ isn't invariant on $TM$ with respect to $\bar{J}.$ Because of this, $\bar{J}(D_2)\subset \Gamma(tr(TM)).$ Thus, it is clear that $\bar{J}(D_1)\neq D_2$ and $\bar{J}(D_1)$ is a distribution on $S(TM).$ Moreover, for $ltr(TM)=L_1 \bot L_2$ and $\xi_1\in\Gamma(D_1),$ $N_1\in\Gamma(L_1),$ since $\bar{g}(\xi_1,N_1)=1,$ then $\bar{J}(L_1)$ is a distribution on $S(TM),$ too.
On the other hand, from (b), $\forall X,Y\in\Gamma(D)$ and $N_1, N_2\in\Gamma(ltr(TM)),$ we get
\begin{align*}
\bar{g}(\bar{R}(X,Y)N_1,N_2)&=\frac{c}{4}\{\bar{g}(Y,N_1)\bar{g}(X,N_2)
-\bar{g}(X,N_1)\bar{g}(Y,N_2) \\
&+\bar{g}(\bar{J}Y,N_1)\bar{g}(\bar{J}X,N_2)
-\bar{g}(\bar{J}X,N_1)\bar{g}(\bar{J}Y,N_2) \\
&+2\bar{g}(X,\bar{J}Y)
\bar{g}(\bar{J}N_1,N_2)\} \\
&=0.
\end{align*}
That is, $\bar{J}(ltr(TM))\cap Rad(TM)\neq 0.$ Thus, $\bar{J}(L_2)$ isn't belong to $Rad(TM)$ or $ltr(TM).$ From this, it is clear that $\bar{J}(D_2)\neq L_2$ and then $\bar{J}D_2\subset S(TM^{\bot}).$ On the other hand, for $\forall \xi_2 \in \Gamma(D_2)$ and $N_2 \in \Gamma(L_2),$ since $\bar{g}(\xi_2,N_2)=1,$ then $\bar{J}(L_2)$ is a distribution on $S(TM^{\bot}),$ too.
Finally, from (c), there exists a non-degenerate distribution $S$ such that $S\bot D$ and for $\forall X,Y\in \Gamma(D)$ and $W_1, W_2 \in\Gamma(S),$ we have
\begin{align*}
\bar{g}(\bar{R}(X,Y)W_1,W_2)&=\frac{c}{4}\{\bar{g}(Y,W_1)  \bar{g}(X,W_2)
-\bar{g}(X,W_1)\bar{g}(Y,W_2) \\
&+\bar{g}(\bar{J}Y,W_1)\bar{g}(\bar{J}X,W_2)
-\bar{g}(\bar{J}X,W_1)\bar{g}(\bar{J}Y,W_2) \\
&+2\bar{g}(X,\bar{J}Y)
\bar{g}(\bar{J}W_1,W_2)\} =0.
\end{align*}
In other words, $S \bot \bar{J}(S).$ Moreover, since $S \bot D$ and $D$ is invariant, we can write
$$
\bar{g}(X,W)=\bar{g}(\bar{J}X,W)=-\bar{g}(X,\bar{J}W)=0,
$$
for $\forall X\in \Gamma(D)$ and $W \in\Gamma(S),$ that is, $\bar{J}(S)$ is orthogonal to $D,$ too. Thus, $S$ and $\bar{J}(S)$ are distributions on $S(TM^{\bot})$ and $S(TM),$ respectively. Thus, $M$ is a STCR lightlike submanifold of $\bar{M}$ and proof is completed. \\

We now investigate the geometry of various distributions defined on $M.$ \\

\noindent{\bf Theorem~3.3.~}{\sl Let $M$ be a screen transversal Cauchy Riemann lightlike submanifold of an indefinite K\"{a}hler manifold $\bar{M}.$ Then
\begin{enumerate}
  \item [(i)] {\sl The distribution $D$ is integrable if and only
  if}
 $$
  h(X, \bar{J}Y)=h(\bar{J}X,Y)\,\,,\,\ \forall X,Y \in \Gamma(D).
 $$
  \item [(ii)] {\sl The distribution $\tilde{D}$ is integrable if and only
  if}
 $$
  A_{\bar{J}Z}V=A_{\bar{J}V}Z \,\,,\,\ \forall Z,V \in \Gamma(\tilde{D}).
 $$
\end{enumerate} }
\noindent{\bf Proof.~} We only prove (i), (ii) is similar. From (\ref{eq:3.12}) and (\ref{eq:3.13}), for $\forall X,Y \in \Gamma(D),$ we
have
\begin{align*}
\omega(\nabla_X Y)=h(X,TY)-Ch(X,Y).
\end{align*}
Hence, since $h$ is symmetric, we obtain
\begin{eqnarray*}
&&\omega([X,Y])=h(X,TY)-h(TX,Y)
\end{eqnarray*}
which proves (i). \\

For the distribution $D,$ we have the following integrability conditions.\\

\noindent{\bf Theorem~3.4.~}{\sl Let $M$ be a STCR lightlike submanifold of an indefinite K\"{a}hler manifold $\bar{M}.$ Then, $D$ is integrable if and only if, for $\forall X,Y\in\Gamma(D),$ the followings are hold:}
\begin{equation*}
h^{s}(X,\bar{J}Y)-h^{s}(Y,\bar{J}X)\in\Gamma (\bar{J}(L_2))
\end{equation*}
and
\begin{equation*}
h^{l}(X,\bar{J}Y)-h^{l}(Y,\bar{J}X)\in\Gamma (L_2).
\end{equation*}
\noindent{\bf Proof.~} We know that $D$ is integrable iff for $\forall X,Y\in\Gamma(D),$ $[X,Y]\in\Gamma(D),$ i.e.,
$$\bar{g}([X,Y],N_{2})=\bar{g}([X,Y],\bar{J}\xi_{1})=\bar{g}([X,Y],\bar{J}W)=0.$$
Thus, $\forall X,Y \in \Gamma(D)$ and $N_{2}\in\Gamma(L_{2}),$ using (\ref{eq:2.2}) and (\ref{eq:2.4}) we have
\begin{align}
\bar{g}([X,Y],N_{2})&= \bar{g}(h^{s}(X,\bar{J}Y)-h^{s}(Y,\bar{J}X),\bar{J}N_2). \label{eq:3.17}
\end{align}
For $\forall X,Y \in\Gamma(D),$ $\xi_{1}\in\Gamma(D_{1})$ and $W\in\Gamma(S),$ using again (\ref{eq:2.2}) and (\ref{eq:2.4}) we obtain
\begin{align}
\bar{g}([X,Y],\bar{J}\xi_{1})&=\bar{g}(h^{l}(Y,\bar{J}X)-h^{l}(X,\bar{J}Y),\xi_1) \label{eq:3.18}
\end{align}
and
\begin{equation}
\bar{g}([X,Y],\bar{J}W)=\bar{g}( h^{s}(Y,\bar{J}X)-h^{s}(X,\bar{J}Y),W).  \label{eq:3.19}
\end{equation}
Thus, from (\ref{eq:3.17})-(\ref{eq:3.19}), the proof is completed. \\

\noindent{\bf Theorem~3.5.~}{\sl Let $M$ be a STCR lightlike submanifold of an indefinite K\"{a}hler manifold $\bar{M}.$ Then, $\tilde{D}$ is integrable if and only if for $\forall X,Y\in\Gamma(\tilde{D}),$ }
\begin{equation*}
A_{\bar{J}X} Y-A_{\bar{J}Y} X\in\Gamma (\tilde{D}).
\end{equation*}
\noindent{\bf Proof.~} $\tilde{D}$ is integrable iff for $\forall X,Y\in\Gamma(\tilde{D}),$ $[X,Y]\in\Gamma(\tilde{D}),$ i.e.,
$$\bar{g}([X,Y],N_{1})=\bar{g}([X,Y],\bar{J}N_{1})=\bar{g}([X,Y],Z)=0.$$
Thus, $\forall X,Y \in \Gamma(\tilde{D}),$ from (\ref{eq:2.1}), (\ref{eq:2.2}), (\ref{eq:2.7}) and (\ref{eq:2.8}) we derive
\begin{align}
\bar{g}([X,Y],N_{1})&= \bar{g}(A_{\bar{J}X} Y-A_{\bar{J}Y} X,\bar{J}N_1). \label{eq:3.20}
\end{align}
For $\forall X,Y \in\Gamma(\tilde{D}),N_{1}\in\Gamma(L_{1}),$ from (\ref{eq:2.2}) and (\ref{eq:2.7}) we obtain
\begin{align}
\bar{g}([X,Y],\bar{J}N_{1})&= \bar{g}(A_{\bar{J}X} Y-A_{\bar{J}Y} X,N_1). \label{eq:3.21}
\end{align}
In a similar way, for $\forall X,Y \in\Gamma(\tilde{D})$ and $Z\in\Gamma(D_0),$ we get
\begin{align}
\bar{g}([X,Y],Z)&=\bar{g}(A_{\bar{J}X} Y-A_{\bar{J}Y} X,\bar{J}Z)=0. \label{eq:3.22}
\end{align}
Then from (\ref{eq:3.20})-(\ref{eq:3.22}) the proof is completed.\\

We now study the geometry of leaves.\\

\noindent{\bf Theorem~3.6.~}{\sl Let $M$ be a STCR lightlike submanifold of an indefinite K\"{a}hler manifold $\bar{M}.$ Then, $D$ defines a totally geodesic foliation on $M$ iff}
$$Bh(X,\bar{J}Y)=0\,,\,\,\ \forall X,Y\in\Gamma(D).$$
\noindent{\bf Proof.~} We assume that $D$ defines a totally geodesic foliation on $M.$ That is, for $\forall X,Y\in\Gamma(D),$ $\nabla_X Y\in\Gamma(D).$ Then, using (\ref{eq:2.4}) and (\ref{eq:2.2}) for $\forall X,Y\in\Gamma(D),$ $\xi_1\in\Gamma(D_{1}),$ $N_2\in\Gamma(L_{2})$ and $W\in\Gamma(S),$ we obtain
\begin{align}
g(\nabla_X Y,\bar{J}\xi_1)&=-\bar{g}(h^{l}(X,\bar{J}Y),\xi_1)=0, \label{eq:3.23} \\
\bar{g}(\nabla_X Y,N_2)&=\bar{g}(h^{s}(X,\bar{J}Y),\bar{J}N_2)=0, \label{eq:3.24} \\
g(\nabla_X Y,\bar{J}W)&=-\bar{g}(h^{s}(X,\bar{J}Y),W)=0. \label{eq:3.25}
\end{align}
Thus, from (\ref{eq:3.23})-(\ref{eq:3.25}) it is easy to see that
$h^{l}(X,\bar{J}Y)$ has no components in $\Gamma(L_1)$ and $h^{s}(X,\bar{J}Y)$ has no components in $\Gamma(L_2 \cup S),$ in other words $Jh(X,Y)$ has no components in $\Gamma(TM)$ and the proof is completed.\\

\noindent{\bf Theorem~3.7.~}{\sl Let $M$ be a STCR lightlike submanifold of an indefinite K\"{a}hler manifold $\bar{M}.$ Then, $\tilde{D}$ defines a totally geodesic foliation on $M$ iff}
$$A_{\omega Y} X\in\Gamma(\tilde{D})\,,\,\,\ \forall X,Y\in\Gamma(\tilde{D}).$$
\noindent{\bf Proof.~} We assume that $\tilde{D}$ defines a totally geodesic foliation on $M.$ That is, for $\forall X,Y \in\Gamma(\tilde{D}),$ $\nabla_X Y\in\Gamma(\tilde{D}).$ Since $Y$ and $\nabla_X Y$ are belong to $\in\Gamma(\tilde{D}),$ then both $TY$ and $T(\nabla_X Y)$ are zero and from (\ref{eq:3.11}) we get
\begin{equation*}
\nabla_X TY-T\nabla_X Y=A_{\omega_L Y} X+A_{\omega_S Y} X+A_{\omega_{\bar{2}} Y}X+Bh(X,Y).
\end{equation*}
From here,
\begin{equation*}
-Bh(X,Y)=A_{\omega_Y}X \in\Gamma(\tilde{D})
\end{equation*}
is obtained. Conversely, if $A_{\omega Y} X\in\Gamma(\tilde{D}),$ for $\forall X,Y\in\Gamma(\tilde{D}),$ then again from (\ref{eq:3.11}), since
$$-A_{\omega_Y}X=T\nabla_X Y+Bh(X,Y),$$
we obtain
$$T\nabla_X Y=0$$
which completes the proof. \\

{\sl As in the Riemannian \cite{Bejancu} and CR-lightlike cases \cite{Sahin-Gunes}, we say that $M$ is a $D$-geodesic (or $\tilde{D}$-geodesic) STCR lightlike submanifold if its second fundamental form $h$ satisfies
\begin{equation}
h(X,Y)=0 \,,\,\,\ \forall X,Y\in\Gamma(D)\,,\,\,\ (or \,\,\ \forall X,Y\in\Gamma(\tilde{D})). \label{eq:3.a}
\end{equation}
It is easy to see that $M$ is a $D$-geodesic (or $\tilde{D}$-geodesic) STCR lightlike submanifold if
\begin{equation}
h^{l}(X,Y)=0 \,,\,\,\ h^{s}(X,Y)=0 \,,\,\,\ \forall X,Y\in\Gamma(D)\,,\,\,\ (or \,\,\ \forall X,Y\in\Gamma(\tilde{D})). \label{eq:3.b}
\end{equation}}

\noindent{\bf Theorem~3.8.~}{\sl Let $M$ be a STCR lightlike submanifold of an indefinite K\"{a}hler manifold $\bar{M}.$ Then, $D$ defines a totally geodesic foliation on $\bar{M}$ iff $M$ is $D$-geodesic.} \\
\noindent{\bf Proof.~} $\Rightarrow):$ We assume that $D$ defines a totally geodesic foliation on $\bar{M}.$ Then, using (\ref{eq:2.4}) for $\forall X,Y \in\Gamma(D),$ $\xi\in\Gamma(Rad(TM))$ and $W\in\Gamma(S(TM^{\bot})),$
\begin{align}
\bar{g}(\bar{\nabla}_X Y,\xi)&=\bar{g}(h^{l}(X,Y),\xi)=0, \label{eq:3.26} \\
\bar{g}(\bar{\nabla}_X Y,W)&=\bar{g}(h^{s}(X,Y),W)=0  \label{eq:3.27}
\end{align}
is obtained. Thus, if $D$ is a totally geodesic foliation on $\bar{M},$ then it is clear that from (\ref{eq:3.26}) and (\ref{eq:3.27}), $h(X,Y)=0$ on $D.$ In other words, since $h(X,Y)=0,$ $\forall X,Y \in\Gamma(D),$ then $M$ is $D$-geodesic. \\
$\Leftarrow ):$ Conversely, let $M$ be $D$-geodesic. Using (\ref{eq:3.a}), (\ref{eq:2.4}) and (\ref{eq:2.2}) we get
\begin{align}
\bar{g}(\bar{\nabla}_X Y,\bar{J}\xi_2)&=-\bar{g}(h^{l}(X,\bar{J}Y),\xi_2)=0, \label{eq:3.28} \\
\bar{g}(\bar{\nabla}_X Y,\bar{J}W)&=-\bar{g}(h^{s}(X,\bar{J}Y),W)=0, \label{eq:3.29}
\end{align}
for $\xi_2\in\Gamma(D_2)\,,\,\ W\in\Gamma(S).$ Thus, from (\ref{eq:3.28}) and (\ref{eq:3.29}) we obtain $\bar{\nabla}_X Y \in\Gamma(D),$ $\forall X,Y \in\Gamma(D)$ and the proof is completed. \\

{\sl We say that $M$ is mixed geodesic STCR lightlike submanifold if its second fundamental form $h$ satisfies
\begin{equation}
h(X,Y)=0 \,,\,\,\ \forall X\in\Gamma(D)\,,\,\ Y\in\Gamma(\tilde{D}). \label{eq:3.30}
\end{equation}
It is easy to see that $M$ is a mixed geodesic STCR lightlike submanifold if
\begin{equation}
h^{l}(X,Y)=0 \,,\,\,\ h^{s}(X,Y)=0 \,,\,\,\ \forall X\in\Gamma(D)\,,\,\ Y\in\Gamma(\tilde{D}). \label{eq:3.31}
\end{equation} }

In the sequel, we find necessary and sufficient conditions for a STCR lightlike submanifold to be mixed geodesic.\\

\noindent{\bf Proposition~3.3~} {\sl Let $M$ be a STCR lightlike submanifold of an indefinite K\"{a}hler manifold $\bar{M}.$ Then, $M$ is mixed geodesic if and only if
\begin{eqnarray*}
A_{\bar{J}Z} X\in\Gamma(D)\,\,\ and \,\,\ \nabla^{t}_X \bar{J}Z\in\Gamma(L_1 \bot D_2 \bot S),
\end{eqnarray*}
for $\forall X\in\Gamma(D)\,,\,\ Z\in\Gamma(\tilde{D}).$ }\\
\noindent{\bf Proof.~} From (\ref{eq:2.4}) and (\ref{eq:2.2}) we obtain
\begin{equation*}
h(X,Z)=-\bar{J}(-A_{\bar{J}Z} X+\nabla^{t}_X \bar{J}Z)-\nabla_X Z, \\
\end{equation*}
for $\forall X\in\Gamma(D)\,,\,\ Z\in\Gamma(\tilde{D}).$ Then, using (\ref{eq:3.8}), (\ref{eq:3.9}) and transversal part we get
\begin{equation*}
h(X,Z)=\omega(A_{\bar{J}Z} X)+C\nabla^{t}_X \bar{J}Z=0 \\
\end{equation*}
which completes the proof. \\

\noindent{\bf Definition~2.~} {\sl A screen transversal Cauchy Riemann lightlike submanifold $M$ of an indefinite K\"{a}hler manifold $\bar{M}$ is called as mixed foliated screen transversal Cauchy Riemann lightlike submanifold if the following conditions are satisfied.
\begin{enumerate}
\item [(1)] $h(X,Z)=0 \,,\,\,\ \forall X\in\Gamma(D)\,,\,\ Z\in\Gamma(\tilde{D}).$
\item [(2)] Distribution $D$ is integrable. \\
\end{enumerate} }

\noindent{\bf Theorem~3.9.~}{\sl Let $M$ be a mixed foliated screen transversal Cauchy Riemann lightlike submanifold of an indefinite K\"{a}hler manifold $\bar{M}.$ Then,
\begin{enumerate}
\item [(i)] Distribution $D$ is parallel iff $Q\nabla^{*}_X Y+Qh^{*}(X,Y)=0,$  $\forall X\in\Gamma(TM)\,,\,\ Y\in\Gamma(D).$
\item [(ii)] Distribution $\tilde{D}$ is parallel iff $P\nabla^{*}_X Z+Ph^{*}(X,Z)=0,$  $\forall X\in\Gamma(TM)\,,\,\ Z\in\Gamma(\tilde{D}).$
\end{enumerate} }
\noindent{\bf Proof.~} \begin{enumerate}
\item [(i)] From (\ref{eq:3.6}) and (\ref{eq:2.11}) we have
\begin{equation*}
\nabla_X Y=P\nabla^{*}_X Y+Ph^{*}(X,Y)+Q\nabla^{*}_X Y+Qh^{*}(X,Y),
\end{equation*}
$\forall X\in\Gamma(TM)\,,\,\ Y\in\Gamma(D),$ which satisfies (i). \\
\item [(ii)] As similar in (i), we get
\begin{equation*}
\nabla_X Z=P\nabla^{*}_X Z+Ph^{*}(X,Z)+Q\nabla^{*}_X Z+Qh^{*}(X,Z),
\end{equation*}
which satisfies (ii).\\
\end{enumerate}
Thus we have the following result.\\

\noindent{\bf Proposition~3.4~} {\sl Let $M$ be a mixed foliated screen transversal Cauchy Riemann lightlike submanifold of an indefinite K\"{a}hler manifold $\bar{M}.$ Then,
\begin{enumerate}
\item [(i)] Distribution $D$ is parallel respect to $\bar{\nabla}$ iff $M$ is $D$-geodesic \\ and $D$ is parallel respect to $\nabla.$
\item [(ii)] Distribution $\tilde{D}$ is parallel respect to $\bar{\nabla}$ iff $M$ is $\tilde{D}$-geodesic \\ and $\tilde{D}$ is parallel respect to $\nabla.$
\end{enumerate} }

\section*{{\bf 4.~ Totally umbilical STCR lightlike submanifolds}}

\setcounter{equation}{0}
\renewcommand{\theequation}{4.\arabic{equation}}
In this section we study totally umbilical {\bf STCR}-lightlike submanifolds, give an example and investigate the existence of such submanifolds. \\

\noindent{\bf Definition~3.~} {\sl \cite{Duggal-Jin} A lightlike submanifold $(M,g)$ of a semi-Riemannian manifold $(\bar{M},\bar{g})$ is totally umbilical in $\bar{M}$ if there is a smooth transversal vector field $H\in\Gamma(ltr(TM))$ on $M,$ called the transversal curvature vector field of $M,$ such that, for all $X,Y\in\Gamma(TM),$
\begin{equation}
h(X,Y)=Hg(X,Y) \label{eq:4.1}.
\end{equation} }
Using (\ref{eq:2.1}) and (\ref{eq:2.3}) it is easy to see that $M$ is totally umbilical if and only if on each coordinate neighborhood $U$ there exist smooth vector fields $H^{l}\in\Gamma(ltr(TM))$ and $H^{s}\in\Gamma(S(TM^{\bot})),$ such that
\begin{align}
h^{l}(X,Y)&=H^{l}g(X,Y) \,,\,\ D^{l}(X,W)=0, \nonumber \\
h^{s}(X,Y)=H^{s}g(X,Y)&\,,\,\ \forall X,Y \in\Gamma(TM)\,,\,\ W\in\Gamma(S(TM^{\bot})). \label{eq:4.2}
\end{align}

The above definition does not depend on the $S(TM)$ and $S(TM^{\bot})$ of $M.$ \\

\noindent{\bf Example~4.~} {\sl Let $\bar{M}=R^{10}_{4}$ be a semi-Euclidean space of signature \\ $(+,+,-,-,+,+,-,-,+,+)$ with respect to the canonical basis
 $$({\partial x_1},{\partial x_2},{\partial x_3},{\partial x_4},{\partial x_5},{\partial x_6},{\partial x_7},{\partial x_8},{\partial x_9},{\partial x_{10}}).$$ \\
Consider a complex structure $\bar{J}$ defined by \\
$$\bar{J}(x_1,x_2,x_3,x_4,x_5,x_6,x_7,x_8,x_9,x_{10})=(-x_2,x_1,-x_4,x_3,-x_6,x_5,-x_8,x_7,-x_{10},x_9).$$
Let $M$ be a submanifold of $(R^{10}_{4},\bar{J})$ given by
\begin{align*}
x_1&=u_1 cosh\alpha \,,\,\,\ x_2=e^{u_2}sinh\alpha+(u_3+\frac{u_4}{2})cosh\alpha \,,\,\,\ x_3=u_1 sinh\alpha ,\\
x_4&=e^{u_2}cosh\alpha+(u_3+\frac{u_4}{2})sinh\alpha \,,\,\,\ x_5=0 \,,\,\,\ x_6=e^{u_2} \,,\,\ x_7=u_1 ,\\
x_8&=u_3-\frac{u_4}{2} \,,\,\,\ x_9=-cosu_5 \,,\,\ x_{10}=sinu_5.
\end{align*}
Then $TM$ is spanned by
\begin{align*}
Z_1&=cosh\alpha{\partial x_1}+sinh\alpha{\partial x_3}+{\partial x_7}\,,\,\,\ Z_2=e^{u_2}\{sinh\alpha{\partial x_2}+cosh\alpha{\partial x_4}+{\partial x_6}\}, \\
Z_3&=cosh\alpha{\partial x_2}+sinh\alpha{\partial x_4}+{\partial x_8}\,,\,\,\ Z_4=\frac{1}{2}\{cosh\alpha{\partial x_2}+sinh\alpha{\partial x_4}-{\partial x_8}\}, \\
Z_5&=sinu_5{\partial x_9}+cosu_5{\partial x_{10}}.
\end{align*}
It is clear that $Rad(TM)=Span\{Z_1,Z_2\}$ and we get by direct calculation that the lightlike transversal bundle $ltr(TM)$ is spanned by
\begin{align*}
N_1&=\frac{1}{2}\{cosh\alpha{\partial x_1}+sinh\alpha{\partial x_3}-{\partial x_7}\}=\bar{J}Z_4 ,\\
N_2&=\frac{1}{2}e^{-u_2}\{-sinh\alpha{\partial x_2}-cosh\alpha{\partial x_4}+{\partial x_6}\}.
\end{align*}
Hence, we have
\begin{align*}
D_1&=Sp\{Z_1\}\,,\,\ D_2=Sp\{Z_2\} \,,\,\ L_1=Sp\{N_1\}\,,\,\ L_2=Sp\{N_2\}\,,\ \\
\bar{J}(D_1)&=Sp\{\bar{J}Z_1=Z_3\}\,,\,\ \bar{J}(L_1)=Sp\{\bar{J}N_1=Z_4\}\,,\ \\
\bar{J}(S)&=Sp\{Z_5\}\,,\,\ S(TM^{\bot})=Sp\{\bar{J}Z_2,\bar{J}N_2,\bar{J}Z_5\}.
\end{align*}
Thus, $M$ be a screen transversal lightlike submanifold of $\bar{M}.$
Moreover, since $\bar{\nabla}_{Z_i} {Z_j}=0,$ $$h^{l}(Z_i,Z_j)=h^{s}(Z_i,Z_j)=0 \,,\,\ 1\leq i\leq 5 \,,\,\ 1\leq j\leq 4 $$
and using Gauss-Weingarten equations
$$h^{l}(Z_5,Z_5)=0 \,,\,\ h^{s}(Z_5,Z_5)=-\bar{J}Z_5$$
is obtained. If we choose $H^{s}=-\bar{J}Z_5,$ then we have $h^{s}(Z_5,Z_5)=H^{s}g(Z_5,Z_5),$ that is, $M$ is a totally umbilical screen transversal lightlike submanifold of $\bar{M}.$ } \\

\noindent{\bf Theorem~4.1.~}{\sl Let $M$ be a totally umbilical STCR lightlike submanifold of an indefinite K\"{a}hler manifold $\bar{M}.$ If distribution $D_0$ is integrable then the induced connection $\nabla$ is a metric connection.} \\
\noindent{\bf Proof.~} Let $M$ be totally umbilical and $D_0$ be integrable. Then, using (\ref{eq:2.2}) we get
\begin{equation*}
\nabla_X \bar{J}Y+h^{l}(X,\bar{J}Y)+h^{s}(X,\bar{J}Y)=\bar{J}\nabla_X Y
+\bar{J}h^{l}(X,Y)+\bar{J}h^{s}(X,Y).
\end{equation*}
Taking lightlike transversal parts of this equation, we have
\begin{equation}
g(X,\bar{J}Y)H^{l}=\omega_L {\nabla_X Y}+g(X,Y)CH^{l}, \label{eq:4.3}
\end{equation}
$\forall X,Y \in\Gamma(D_0).$ Hence from (\ref{eq:4.3})
\begin{equation}
g(X,\bar{J}Y)H^{l}-g(Y,\bar{J}X)H^{l}+\omega_L([X,Y])=0 \label{eq:4.4}
\end{equation}
is obtained. Since $D_0$ is integrable, for $\forall X, Y \in\Gamma(D_0),$ then $\nabla_X Y \in\Gamma(D_0).$ If we choose $X=\bar{J}Y \in\Gamma(D_0),$ from (\ref{eq:3.1}) we get $$2g(Y,Y)H^{l}=0.$$ We know that since $D_0$ is non-degenerate, then $H^{l}=0.$ That is, from (\ref{eq:4.2}) $h^{l}=0.$ Hence, from (\ref{eq:2.16}) the proof is completed. \\

\noindent{\bf Theorem~4.2.~}{\sl Let $M$ be a totally umbilical STCR lightlike submanifold of an indefinite K\"{a}hler manifold $\bar{M}.$ Then, $H^{s}\in\Gamma(\bar{J}(D_2) \cup S).$} \\
\noindent{\bf Proof.~} Let $M$ be totally umbilical STCR lightlike submanifold. Then for $\forall X,Y \in\Gamma(D_0),$ if we consider equation (\ref{eq:3.13}), we get
\begin{equation}
h^{s}(X,\bar{J}Y)=Ch^{s}(X,Y)+\omega_S (\nabla_X Y)+\omega_{\bar{2}} (\nabla_X Y). \label{eq:4.5}
\end{equation}
If we choose $X=Y\in\Gamma(D_0)$ in (\ref{eq:4.5}), we obtain $CH^{s}=0$ which completes the proof. \\

\noindent{\bf Theorem~4.3.~}{\sl Let $M$ be a totally umbilical STCR lightlike submanifold of an indefinite K\"{a}hler manifold $\bar{M}.$ Then, one of the followings is hold: }
\begin{enumerate}
  \item [(a)] $M$ is totally geodesic.
  \item [(b)] $h^{s}=0$ or $dim(S)=1$ and $D_0$ is not integrable.
\end{enumerate}
\noindent{\bf Proof.~} If $D_0\neq \{0\}$ and integrable, from Theorem 4.1 we obtain that $h^{l}=h^{s}=0$ which is case (a). Now suppose that $D_0$ is not integrable. Then using (\ref{eq:2.2}), (\ref{eq:2.8}) and taking tangential part we obtain
\begin{equation}
-A_{\bar{J}W} Z=T\nabla_Z W+Bh(Z,W), \label{eq:4.6}
\end{equation}
$Z,W\in\Gamma(\bar{J}S).$ On the other hand, if we take $X=Y=Z$ and $W=\bar{J}W$ in (\ref{eq:2.9}) and use (\ref{eq:4.6}),
\begin{equation}
\bar{g}(h^{s}(Z,Z)\bar{J}W)=-\bar{g}(h(Z,W),\bar{J}Z), \label{eq:4.7}
\end{equation}
$\forall Z,W\in\Gamma(\bar{J}(S)),$ is obtained.
Since $M$ is totally umbilical, from (\ref{eq:4.7}) we have
\begin{equation}
g(Z,Z)\bar{g}(H^{s},\bar{J}W)=-\bar{g}(Z,W)\bar{g}(H^{s},\bar{J}Z). \label{eq:4.8}
\end{equation}
Interchanging the role of $Z$ and $W$ in this equation we get
\begin{equation}
\bar{g}(H^{s},\bar{J}Z)=\frac{\bar{g}(Z,W)^{2}}{g(Z,Z)g(W,W)} \bar{g}(H^{s},\bar{J}Z). \label{eq:4.9}
\end{equation}
Since $\bar{J}(S)$ is non-degenerate, choosing non-null vector fields $Z$ and $W,$ we conclude that either $H^{s}=0$ or $Z$ and $W$ are linearly dependent. This proves (b). Thus, the proof is completed.  \\

\noindent{\bf Theorem~4.4.~}{\sl There exist no totally umbilical proper STCR lightlike submanifold of an indefinite complex space form $\bar{M}(c),$ $c\neq 0.$} \\
\noindent{\bf Proof.~} Let assume that $M$ be a totally umbilial proper STCR lightlike submanifold of an indefinite complex space form $\bar{M}(c),$ $c\neq 0.$ Then from (\ref{eq:2.3}) and (\ref{eq:2.17}) we obtain
\begin{equation}
\bar{R}(X,\bar{J}X)Z=-\frac{c}{2}g(X,X)\bar{J}Z \label{eq:4.10}
\end{equation}
and
\begin{align}
\bar{R}(X,\bar{J}X)Z&=(\nabla_X h^{s})(\bar{J}X,Z)-(\nabla_{\bar{J}X} h^{s})(X,Z) \nonumber \\
&=Xh^{s}(\bar{J}X,Z)-h^{s}(\nabla_X \bar{J}X,Z)-h^{s}(\bar{J}X,\nabla_X Z) \nonumber \\
&-\bar{J}Xh^{s}(X,Z)+h^{s}(\nabla_{\bar{J}X} X,Z)+h^{s}(X,\nabla_{\bar{J}X} Z) \nonumber \\
&=Xg(\bar{J}X,Z)H^{s}-g(\nabla_X \bar{J}X,Z)H^{s}-g(\bar{J}X,\nabla_X Z)H^{s} \nonumber \\
&-\bar{J}Xg(X,Z)H^{s}+g(\nabla_{\bar{J}X} X,Z)H^{s}+g(X,\nabla_{\bar{J}X} Z)H^{s} \nonumber \\
&=0, \label{eq:4.11}
\end{align}
for $\forall X \in\Gamma(D_0) \,,\,\ Z\in\Gamma(\bar{J}S),$  respectively.
Thus, from (\ref{eq:4.10}) and (\ref{eq:4.11}) we have $c=0$ which is a contraction and the proof is completed. \\

\section*{{\bf 5.~ Minimal STCR lightlike submanifolds}}

\setcounter{equation}{0}
\renewcommand{\theequation}{5.\arabic{equation}}

A general notion of minimal lightlike submanifold $M$ of a semi-Riemannian manifold $\bar{M}$ has been introduced by Bejan-Duggal in \cite{Bejan-Duggal} as follows: \\

\noindent{\bf Definition~4.~} {\sl We say that a lightlike submanifold $(M,g,S(TM))$ isometrically immersed in a semi-Riemannian manifold $(\bar{M},\bar{g})$ is minimal if:
\begin{enumerate}
\item [(i)] $h^{s}=0$ on $Rad(TM)$ and
\item [(ii)] $traceh=0,$ where trace is written with respect to $g$ restricted to $S(TM).$
\end{enumerate} }

In Case 2, condition (i) is trivial. It has been shown in \cite{Bejan-Duggal} that the above definition is independent of $S(TM)$ and $S(TM^{\bot}),$ but it depends on $tr(TM).$ \\

As in the semi-Riemannian case, any lightlike totally geodesic $M$ is minimal. Thus, it follows from Corollary 2.5 in \cite[page 167]{Duggal-Bejancu} that any totally lightlike $M$ (Case 4) is minimal. If $M$ is totally umbilical proper STCR lightlike submanifold of an indefinite K\"{a}hler manifold $\bar{M}$ with the distribution $D_0$ integrable, then it follows from Theorem 4.2 that $M$ is minimal. \\
Minimal lightlike submanifolds are investigated in detail in \cite{Duggal-Sahin}. \\

\noindent{\bf Example~5.~} {\sl Let $M$ be a submanifold of $\bar{M}=R^{14}_{4}$ and given by
\begin{align*}
x_1&=u_1 \sinh\alpha+u_2 \cosh\alpha \,,\,\ x_2=(u_3+\frac{u_4}{2})\sinh\alpha \,,\,\ x_3=u_1 \,,\,\ x_4=(u_3-\frac{u_4}{2})\,,\ \\
x_5&=u_2 \,,\,\ x_6=0 \,,\,\ x_7=\cos u_5 \cosh u_6 \,,\,\ x_8=\sin u_5 \sinh u_6\,,\ \\
x_9&=u_1 \cosh\alpha+u_2 \sinh\alpha \,,\,\ x_{10}=(u_3+\frac{u_4}{2})\cosh\alpha \,,\,\ x_{11}=-\sin u_7\cosh u_8 \,,\ \\
x_{12}&=\cos u_7 \sinh u_8 \,,\,\ x_{13}=-(u_7+u_8) \,,\,\ x_{14}=u_7-u_8.
\end{align*}
Then $TM$ is spanned by
\begin{align*}
Z_1&=\sinh\alpha \partial x_1+{\partial x_3}+\cosh\alpha {\partial x_9}\,,\,\,\ Z_2=\cosh\alpha \partial x_1+{\partial x_5}+\sinh\alpha {\partial x_9}, \\
Z_3&=\sinh\alpha\partial x_2+{\partial x_4}+\cosh\alpha {\partial x_{10}} \,,\,\,\
Z_4=\frac{1}{2}\{\sinh\alpha\partial x_2+{\partial x_4}+\cosh\alpha {\partial x_{10}}\}, \\
Z_5&=-\sin u_5 \cosh u_6 {\partial x_7}+ \cos u_5 \sinh u_6 {\partial x_8}, \\
Z_6&=\cos u_5\sinh u_6 {\partial x_7}+\sin u_5\cosh u_6 {\partial x_8}, \\
Z_7&=-\cos u_7\cosh u_8 \frac{\partial}{\partial x_{11}}-\sin u_7\sinh u_8 {\partial x_{12}}-{\partial x_{13}}+{\partial x_{14}}, \\
Z_8&=-\sin u_7\sinh u_8 {\partial x_{11}}+\cos u_7\cosh u_8 {\partial x_{12}}-{\partial x_{13}}-{\partial x_{14}}.
\end{align*}
where
$Rad(TM)=Sp\{{Z_1,Z_2}\} \,,\,\,\ \bar{J}(D_1)=Sp\{{\bar{J}Z_1=Z_3}\} \,,\,\,\ D_0=Sp\{{Z_5,Z_6}\}$
and it is easy to see that
\begin{align*}
ltr(TM)=Sp\{&N_1=\frac{1}{2}\{{\sinh\alpha \partial x_1}-{\partial x_3}+{\cosh\alpha\ \partial x_9}\},\\
&N_2=\frac{1}{2}\{-\cosh\alpha \partial x_1+{\partial x_5}-\sinh\alpha {\partial x_9}\}\}, \\
\bar{J}N_1=Z_4 \,,\,\,\ S(T&M^{\bot})=Sp\{\bar{J}Z_2,\bar{J}N_2,\bar{J}Z_7,\bar{J}Z_8\}, \\
&S=Sp\{\bar{J}Z_7,\bar{J}Z_8\}.
\end{align*}
Hence, $M$ is a STCR lightlike submanifold of $R^{14}_4.$ \\
On the other hand, by direct calculation we obtain
$$\bar{\nabla}_{Z_{i}} {Z_{j}}=0 \,,\,\,\ 1\leq i\leq4 \,,\,\,\ 1\leq j\leq8$$
and
\begin{align*}
h(Z_5,Z_5)&=0 \,,\,\,\ h(Z_6,Z_6)=0, \\
h^{l}(Z_7,Z_7)&=0 \,,\,\,\ h^{l}(Z_8,Z_8)=0, \\
h^{s}(Z_7,Z_7)&=\frac{(\sinh u_8\cosh u_8)\bar{J}Z_7+(-\sin u_7\cos u_7)
\bar{J}Z_8}{\sin^{2}u_7\sinh^{2}u_8+\cos^{2}u_7\cosh^{2}u_8+2}, \\
h^{s}(Z_8,Z_8)&=\frac{(-\sinh u_8\cosh u_8)\bar{J}Z_7+(\sin u_7\cos u_7)
\bar{J}Z_8}{\sin^{2}u_7\sinh^{2}u_8+\cos^{2}u_7\cosh^{2}u_8+2}.
\end{align*}
Hence, it is clear that $M$ is not totally geodesic and, but it is a minimal STCR lightlike submanifold of $\bar{M}=R^{14}_{4}.$ } \\

We now give characterizations for $M$ to be minimal.\\

\noindent{\bf Theorem~5.1.~}{\sl A totally umbilical STCR lightlike submanifold $M$ is minimal if and only if for $W_j\in\Gamma(S(TM^{\bot})),$
\begin{equation*}
trace A_{W_j}\mid_{D_0 \bot \bar{J}(S)}=trace A_{\xi_k}^{*}\mid_{D_0 \bot \bar{J}(S)}=0,
\end{equation*}
where $dim(TM)=m,$ $dim(Rad(TM))=r,$ $dim(tr(TM))=n.$ } \\
\noindent{\bf Proof.~} $M$ is minimal iff $h^{s}=0$ on $Rad(TM)$ and $traceh=0$ on $S(TM).$ That is,
\begin{align*}
traceh\mid_{S(TM)}&=traceh\mid_{D_0} +traceh\mid_{\bar{J}(S)} +traceh\mid_{\bar{J}(D_1)}+traceh\mid_{\bar{J}(L_1)} \\
&=\sum^{a}_{i=1} h(Z_i,Z_i)+\sum^{b}_{j=1} h(\bar{J}W_j,\bar{J}W_j) \\
&+\sum^{r}_{k=1} h(\bar{J}\xi_k,\bar{J}\xi_k)+\sum^{r}_{k=1} h(\bar{J}N_k,\bar{J}N_k),
\end{align*}
where $dim(D_0)=a,$ $dim(\bar{J}(S))=b.$ \\
On the other hand, since $M$ is totally umbilical then from (\ref{eq:4.2}) we get
\begin{align*}
traceh\mid_S(TM)&=traceh\mid_{D_0} +traceh\mid_{\bar{J}(S)} \\
&=\sum^{a}_{i=1} \varepsilon_i \{h^{l}(e_i,e_i)+h^{s}(e_i,e_i)\}+\sum^{b}_{j=1}\varepsilon_j \{h^{l}(e_j,e_j)+h^{s}(e_j,e_j)\} \\
&=\sum^{a}_{i=1} \varepsilon_i [\frac{1}{r}\sum^{r}_{k=1} \bar{g}(h^{l}(e_i,e_i),\xi_k)N_k+\frac{1}{n-r}\sum^{n-r}_{j=1}\bar{g}(h^{s}(e_i,e_i),W_j)W_j] \\
&+\sum^{b}_{j=1} \varepsilon_j [\frac{1}{r}\sum^{r}_{k=1} \bar{g}(h^{l}(e_j,e_j),\xi_k)N_k+\frac{1}{n-r}\sum^{n-r}_{j=1}\bar{g}(h^{s}(e_j,e_j),W_j)W_j].
\end{align*}
Besides, if we consider (\ref{eq:2.9}) and (\ref{eq:2.13}), we obtain
\begin{align*}
traceh\mid_S(TM)&=\sum^{a}_{i=1} \varepsilon_i [\frac{1}{r}\sum^{r}_{k=1} g(A_{\xi_k}^{*} e_i,e_i)N_k]+\sum^{b}_{j=1} \varepsilon_j [\frac{1}{r}\sum^{r}_{k=1}g(A_{\xi_k}^{*} e_j,e_j)N_k] \\
&+\sum^{a}_{i=1} \varepsilon_i [\frac{1}{n-r}\sum^{n-r}_{j=1}g(A_{W_j} e_i,e_i)W_j] \\
&+\sum^{b}_{j=1} \varepsilon_j[\frac{1}{n-r}\sum^{n-r}_{j=1}g(A_{W_j} e_j,e_j)W_j] \\
&=0
\end{align*}
which completes the proof. \\

\noindent{\bf Theorem~5.2.~}{\sl Let $M$ be a STCR lightlike submanifold of an indefinite K\"{a}hler manifold $\bar{M}.$ Then the distribution $D_0$ is minimal if and only if for $\forall X \in\Gamma(D_0)$ and $N \in\Gamma(L_1),$
$$A_N \bar{J}X+\bar{J}A_N X \,\,\ has \,\,\ no \,\,\ components \,\,\ in \,\,\ \Gamma(D_0).$$ }
\noindent{\bf Proof.~} From definition, it is clear that $D_0$ is minimal if and only if
\begin{align*}
g(\nabla_X X+\nabla_{\bar{J}X} \bar{J}X,\bar{J}\xi)&=0, \\
g(\nabla_X X+\nabla_{\bar{J}X} \bar{J}X,\bar{J}W)&=0, \\
g(\nabla_X X+\nabla_{\bar{J}X} \bar{J}X,\bar{J}N^{'})&=0,
\end{align*}
 $\forall X \in\Gamma(D_0)$, $\xi\in\Gamma(Rad(TM)),$ $W\in\Gamma(S)$ and $N^{'}\in\Gamma(L_2).$ From (\ref{eq:2.1}) we can write
\begin{align*}
g(\nabla_X X,\bar{J}\xi)&=-g(\bar{J}X,A^{*} _\xi X), \\
g(\nabla_{\bar{J}X} \bar{J}X,\bar{J}\xi)&=g(X,A^{*} _\xi \bar{J}X).
\end{align*}
On the other hand, the shape operator is symmetric on $S(TM).$ Thus, from (\ref{eq:2.6})-(\ref{eq:2.8}), for $\forall X \in\Gamma(D_0),$ $\xi\in\Gamma(Rad(TM)),$ $W\in\Gamma(S)$ and $N^{'}\in\Gamma(L_2),$ we have
\begin{align*}
g(\nabla_X X+g(\nabla_{\bar{J}X} \bar{J}X,\bar{J}\xi)&=g(A^{*} _\xi X, \bar{J}X)-g(\bar{J}X,A^{*} _\xi X)=0.
\end{align*}
In a similar way, we get
$$g(\nabla_X X+\nabla_{\bar{J}X} \bar{J}X,\bar{J}W)=0.$$
Similarly, for $\forall X \in\Gamma(D_0) \,\,\ and \,\,\ N\in\Gamma(L_1),$ we obtain
\begin{equation}
g(\nabla_X X+\nabla_{\bar{J}X} \bar{J}X,\bar{J}N)=g(X,A_N \bar{J}X+\bar{J}A_N X). \label{eq:5.1}
\end{equation}
Hence, proof is completed.\\

\end{document}